\documentclass[12pt]{article}
\usepackage{amsmath,amsthm}
\usepackage{amssymb,latexsym}
\usepackage{enumerate}
\usepackage[english]{babel}
\usepackage{graphicx}

\def\Q{{\mathbb Q}}
\def\Z{{\mathbb Z}}

\def\O{{\cal O}}

\newtheorem{lemma}{Lemma}
\newtheorem{theorem}[lemma]{Theorem}

\newtheorem{proposition}[lemma]{Proposition}

\title{
Power integral bases in a family of sextic fields with quadratic subfields\\
}
\author{
Istv\'{a}n Ga\'{a}l\thanks{
        Research supported in part by K115479 from the
        Hungarian National Foundation for Scientific Research
                         },\; 
and L\'aszl\'o Remete
\\ \\
University of Debrecen, Mathematical Institute \\
H--4010 Debrecen Pf.12., Hungary \\
e--mail: gaal.istvan@unideb.hu, remetel42@gmail.com
}

\begin{document}

\maketitle
\thispagestyle{empty}

\renewcommand{\thefootnote}{}



\renewcommand{\thefootnote}{\arabic{footnote}}
\setcounter{footnote}{0}

\begin{abstract}
Let $M=\Q(i\sqrt{d})$ be any imaginary quadratic field with a positive
square-free $d$. Consider the polynomial
\[
f(x)=x^3-ax^2-(a+3)x-1,
\]
with a parameter $a\in\Z$. 
Let $K=M(\alpha)$, where $\alpha$ is a root of $f$. 
This is an infinite parametric family of sextic fields depending
on two parameters, $a$ and $d$.
Applying relative Thue equations
we determine the relative power integral bases of these sextic fields
over their quadratic subfields. 
Using these results we also determine generators of (absolute) 
power integral bases of the sextic fields. 
\end{abstract}

\section{Introduction}

Monogenity is a classical topic of algebraic number theory c.f. \cite{nark}, \cite{book}. 
Let $M\subset K$ be algebraic number fields
with $[K:M]=n$, denote by $\Z_M$ and $\Z_K$
the rings of integers of $M$ and $K$, respectively. 
An order $\O$ of $\Z_K$ is {\it monogene over} $M$
if $\O=\Z_M[\vartheta]$ with some $\vartheta\in\O$. 
In this case $(1,\vartheta,\ldots,\vartheta^{n-1})$ is
a {\it power integral basis} of $\O$ {\it over} $M$. If $\O=\Z_K$ then we call $K$
{\it monogene over} $M$. In the absolute case if $M=\Q$ then we call $K$ {\it monogene} and 
$(1,\vartheta,\ldots,\vartheta^{n-1})$ a {\it power integral basis} of $K$.

In the relative case $\alpha,\beta\in\O$ are called {\it equivalent} if 
$\alpha=a+\varepsilon \beta$ with an $a\in\Z_M$ and with a unit $\varepsilon\in\Z_M$.
In the absolute case (if $M=\Q$) this simplifies to $\alpha=a\pm \beta$ with
an $a\in\Z$. If $\alpha$ is equivalent to $\beta$ then $\alpha$ generates a 
power integral basis if and only if $\beta$ does: we determine generators of
power integral bases only up to equivalence. 

There are several algorithms for deciding monogenity and determining generators of
power integral bases. These procedures heavily
depend on the degree and other properties of 
number fields, see \cite{book}. 

We also considered monogenity in infinite parametric families of number fields
in the absolute and relative cases, as well, see
I.Ga\'al and G.Lettl \cite{gl},
I.Ga\'al  and T.Szab\'o \cite{gsz3}, \cite{gsz4}.

Recently I.Ga\'al, L.Remete and T.Szab\'o \cite{grsz2} studied the
relation of monogenity and relative monogenity which was already
applied in a family of octic fields by I.Ga\'al and L.Remete \cite{gr}.
In the present paper we utilize similar tools to describe power integral bases
in a well known infinite parametric family of sextic fields.

\section{A parametric family of sextic fields}

Throughout this paper $d$ will be a positive square-free integer. 
Set $M=\Q(i\sqrt{d})$ with
ring of integers $\Z_M$.
Let $\omega=i\sqrt{d}$ if $-d\equiv 2,3 \ (\bmod \ 4)$ and $\omega=(1+i\sqrt{d})/2$ if 
$-d\equiv 1 \ (\bmod \ 4)$.

Consider the polynomial
\begin{equation}
f(x)=x^3-ax^2-(a+3)x-1,
\label{f}
\end{equation}
with a parameter $a\in\Z$. 
These well-known polynomials correspond to the simplest cubic
fields of D.Shanks \cite{shanks} studied by several authors.
The discriminant of the polynomial is
\begin{equation}
D(f)=(a^2+3a+9)^2.
\label{df}
\end{equation}
Let $\alpha$ be a root of $f$. Remark that $\alpha$ depends on the
parameter $a$. 
Our purpose is to determine generators of power integral bases in the 
sextic fields $K=M(\alpha)$, more exactly in the order
\begin{equation}
\O=\Z[1,\alpha,\alpha^2,\omega,\omega\alpha,\omega\alpha^2]
\label{basis}
\end{equation}
of $K$.
The field $K$ is just the composite of $M=\Q(i\sqrt{d})$ and $L=\Q(\alpha)$.
As it is well known (see \cite{nark}) if $\{1,\alpha,\alpha^2\}$ is an
integral basis in $L$ and the discriminants of $L$ and $M$
are coprime, then we have $\O=\Z_K$.
According to \cite{was}, if $a^2+3a+9$ is square-free, then
$\{1,\alpha,\alpha^2\}$ is indeed an integral basis in $L$.

\section{Simplest cubics over imaginary quadratic fields}

We shall use a result of C.Heuberger \cite{heuberger}
on the solutions of 
the parametric family of relative Thue equations
corresponding to the simplest cubic fields.
Consider the equation 
\begin{equation}
Y_1^3-aY_1Y_2^2-(a+3)Y_1Y_2^2-Y_2^3=\varepsilon \;\;\; (Y_1,Y_2\in\Z_M)
\label{rrr1}
\end{equation}
where $\varepsilon$ is a unit in $M$.
This parametric family of relative Thue equation was studied 
by C.Heuberger, A.Peth\H o and R.F.Tichy \cite{hhh} who gave the solutions
for large parameters and by C.Heuberger \cite{heuberger}, who gave the
solutions for all parameters (even for quadratic integer parameters). 
His result was extended by 
P.Kirschenhofer, C.M.Lampl and J. Thuswaldner \cite{kirsch} 
involving also a wider class of rights hand sides.
Here we use a special case of Heuberger's result (see also \cite{honline}):

\begin{lemma}
Let $\omega_3=(1+i\sqrt{3})/2$. 
Up to sign equation (\ref{rrr1}) has 
the following solutions independent of "$a$":
$(Y_1,Y_2)=(1,0),(0,1),(1,-1),(i,-i),(i,0),(0,i),\\
(\omega_3,-\omega_3),(0,\omega_3),(\omega_3,0),(1-\omega_3,0),(0,1-\omega_3),
(1-\omega_3,-1+\omega_3)$\\
and the following solutions depending on "$a$":\\
$(a,Y_1,Y_2)=(-3,9,-2),(-3,7,-9),(-3,2,7),(-1,3,-1),(-1,2,-3),\\
(-1,1,2),
(0,9,-4),(0,5,-9),(0,4,5),(0,2,-1),(0,1,1),(0,1,-2),
(1,9,-5),\\(1,5,4),(1,4,-9),(1,2,-1),(1,1,1),(1,1,-2),
(2,3,-2),(2,2,1),(2,1,-3),\\
(4,9,-7),(4,7,2),(4,2,-9)$.  
\label{lemma1}
\end{lemma}

\section{Results}

Using the results of I.Ga\'al, L.Remete and T.Szab\'o \cite{grsz2}
we show that the generators of relative power integral bases
of $\O$ over $M$ can be obtained from the solutions of the above
relative Thue equation. As a consequence of Lemma \ref{lemma1} we have
the following proposition which will be proved in Section \ref{sec6}:

\begin{proposition}
Up to equivalence all generators of relative power integral bases
of $\O$ over $M$ are of the form $X_1\alpha+X_2\alpha^2$ where 
$X_1=Y_1-aY_2,X_2=Y_2$, and $(Y_1,Y_2)$ is a solution of 
equation (\ref{rrr1}).
\label{t1}
\end{proposition}

\noindent
Our main result is on the generators of (absolute) power integral bases
of $\O$:

\begin{theorem}
The order $\O$ is only monogene for $d=1$ for the following values of "$a$",
otherwise $\O$ is not monogene.
For $d=1$ these values of "$a$" together with the coordinates 
(listed up to sign)
$y\in\Z,X_1,X_2\in\Z_M$ of the
generators 
\[
y\omega+X_1\alpha+X_2\alpha^2
\]
of power integral bases of $\O$ are given by \\
$(a,y,X_1,X_2)=
(-3,1,i+ai,-i),(-2,1,i+ai,-i),
(-1,1,i+ai,-i),\\(0,1,i+ai,-i),
(-3,2,i+ai,-i),(-2,2,i+ai,-i),
(-1,2,i+ai,-i),\\(0,2,i+ai,-i),
(-3,1,i,0),(-2,1,i,0),
(-1,1,i,0),(0,1,i,0),\\
(-3,0,i,0),(-2,0,i,0),
(-1,0,i,0),(0,0,i,0,),
(-2,0,-ia,i),\\(-1,-1,-ia,i),
(0,-2,-ia,i),(-3,1,-ia,i),
(0,-3,-ia,i),\\(-3,0,-ia,i),
(-1,-2,-ia,i),(-2,-1,-ia,i).$
\label{t2}
\end{theorem}

\section{Absolute and relative monogenity}

The discriminant of the basis (\ref{basis}) of $\O$ is
\begin{equation}
D_{\O}=D(f)^2\cdot D_M^3.
\label{ddd}
\end{equation}
Using the method of \cite{grsz2} we shall first determine 
generators of relative power integral bases of $\O$ over $M$ and then
generators of (absolute) power integral bases of $\O$.

Denote by $\omega^{(1)}$ and $\omega^{(2)}$ the conjugates of $\omega\in M$.
Let $\alpha^{(i)}$ ($i=1,2,3$) be the roots of $F$. Any $\vartheta\in\O$
can be written in the form
\begin{equation}
\vartheta=x_0+x_1\alpha+x_2\alpha^2+y_0\omega+y_1\omega\alpha+y_2\omega\alpha^2
\label{koord}
\end{equation}
with $x_0,x_1,x_2,y_0,y_1,y_2\in\Z$.
Set
\begin{equation}
\vartheta^{(i,j)}=x_0+x_1\alpha^{(j)}
+x_2\left(\alpha^{(j)}\right)^2+y_0\omega^{(i)}+y_1\omega^{(i)}\alpha^{(j)}
+y_2\omega^{(i)}\left(\alpha^{(j)}\right)^2
\label{tetakonj}
\end{equation}
for $1\leq i\leq 2,\; 1\leq j\leq 3$. Let
\begin{equation}
I_{\O/M}(\vartheta)=(\O^+ : \Z_M[\vartheta]^+)
=\frac{1}{|(D(F))|}
\cdot
\prod_{i=1}^2\;\;\prod_{1\leq j_1<j_2\leq 3}
 \left|\vartheta^{(i,j_1)}-\vartheta^{(i,j_2)}\right|
\label{relind}
\end{equation}
and 
\begin{equation}
J(\vartheta)=(\Z_M[\vartheta]^+: \Z[\vartheta]^+)
=\frac{1}{(\sqrt{|D_M|})^3}
\cdot
\prod_{j_1=1}^3 \prod_{j_2=1}^3
\left|\vartheta^{(1,j_1)}-\vartheta^{(2,j_2)}\right|.
\label{jjj}
\end{equation}
In our case Proposition 1 of I.Ga\'al, L.Remete and T.Szab\'o \cite{grsz2}
gets the following from:
\begin{lemma}
We have
\begin{equation}
I(\vartheta)=I_{\O/M}(\vartheta) \cdot J(\vartheta)
\label{index}
\end{equation}
where the first factor is the relative index of $\vartheta$ over $M$.
\label{ppp}
\end{lemma}
Obviously, $\vartheta$ generates a power integral basis in $\O$ if and only
if its index (\ref{index}) is equal to 1, that is both the relative index 
(\ref{relind}) and the factor (\ref{jjj}) has to be equal to 1.

\section{Calculating generators of relative power integral bases}
\label{sec6}

{\bf Proof of Proposition \ref{t1}}\\
First we calculate those $\vartheta\in\O$ for which the relative index 
(\ref{relind}) is 1. These elements generate a relative power integral
basis of $\O$ over $M$.
As we have shown in I.Ga\'al \cite{relcubic}, calculating generators of 
relative power integral bases in cubic relative extensions
leads to cubic relative Thue equations. 

Using the notation (\ref{tetakonj})
for $i=1,2$, $j_1<j_2$ we have
\[
\vartheta^{(i,j_1)}-\vartheta^{(i,j_2)}=
(\alpha^{(j_1)}-\alpha^{(j_2)})
\left(  X_1^{(i)} +(\alpha^{(j_1)}+\alpha^{(j_2)}) X_2^{(i)} \right)
\]
with quadratic integers $X_1=x_1+\omega y_1$, $X_2=x_2+\omega y_2$ in $M$.
We have 
$\alpha^{(1)}+\alpha^{(2)}+\alpha^{(3)}=a$, that is if
$j=\{1,2,3\}\setminus \{j_1,j_2\}$ then 
\[
X_1^{(i)} +(\alpha^{(j_1)}+\alpha^{(j_2)}) X_2^{(i)} 
=X_1^{(i)} +(a -\alpha^{(j)})  X_2^{(i)} 
\]
\[
=(X_1^{(i)}+a X_2^{(i)})-\alpha^{(j)} X_2^{(i)} =
Y_1^{(i)}-\alpha^{(j)} Y_2^{(i)}
\]
with 
\begin{equation}
Y_1=X_1+aX_2, \;\; Y_2=X_2. 
\label{xy}
\end{equation}
We have
\[
\prod_{i=1}^2\;\;\prod_{1\leq j_1<j_2\leq 3}
 \left|\vartheta^{(i,j_1)}-\vartheta^{(i,j_2)}\right|
\]
\[
 =
\left(\prod_{i=1}^2\;\;\prod_{1\leq j_1<j_2\leq 3}|\alpha^{(j_1)}-\alpha^{(j_2)}|\right)
\cdot \left| N_{M/Q}(N_{K/M}(Y_1-\alpha Y_2))\right|.
\]
Therefore by (\ref{relind}) the equation
\[
I_{\O/M}(\vartheta)=1
\]
is just equivalent to
\begin{equation}
N_{M/Q}(N_{K/M}(Y_1-\alpha Y_2))=\pm 1.
\label{rrr}
\end{equation}
This is just the relative Thue equation (\ref{rrr1}) corresponding to the polynomial $f$.

If $(Y_1,Y_2)$ is a solution of the relative Thue equation, 
then by (\ref{xy}) we calculate
\[
X_1=Y_1-aY_2,\; X_2=Y_2.
\]
Using this transformation, Lemma \ref{lemma1} implies Proposition \ref{t1}. $\Box$

\section{Calculating generators of power integral bases}

{\bf Proof of Theorem \ref{t2}}\\
Generators of power
integral bases are determined only up to sign, and translation by elements in $\Z$. In view of Lemma \ref{ppp}, given $(X_1,X_2)$ we have to determine 
$A=x_0+\omega y_0$ and the unit $\varepsilon$ in $M$ so that for 
\begin{equation}
\vartheta=x_0+\omega y_0+\varepsilon (X_1\alpha+X_2\alpha^2)
\label{ttee}
\end{equation}
we have
\begin{equation}
J(\vartheta)=1.
\label{j1}
\end{equation}
The index is translation invariant, indeed $J(\vartheta)$ does not depend on $x_0$.
Therefore we only have to determine $y_0$. This seems to be easy, however
several solutions $(X_1,X_2)$ are independent from the parameter $a$ and from
$d$. This means that in several cases equation (\ref{j1}) has three variables
and degree 9.

\vspace{0.5cm}
\noindent
{\bf I.} Consider first the solutions $(Y_1,Y_2)=(1,-1),(1,0),(0,1)$
which are independent from $d$ and $a$.
Let $X_1=Y_1-aY_2$, $X_2=Y_2$.

\noindent
{\bf I.1.} Let $-d\equiv 2,3 \ (\bmod \ 4), -d\neq -1$.
Taking any of $(Y_1,Y_2)=(1,-1),(1,0),(0,1)$ we have the same arguments.
$J(\vartheta)$ is divisible
by $y_0^3$, therefore $y_0=\pm 1$. The remaining factor $J_1$ of 
$J(\vartheta)$ is quadratic in $y_0$, hence we may substitute $y_0=1$
to get $J_2$. It is easily seen that $J_2 \equiv 1 \ (\bmod \ 4)$
both for even and odd values of $a$. Therefore $J_2-1=0$ must be satisfied
for a generator of a power integral basis.
Set $K=a^2+3a+9$ then
\[
J_2-1=K^2(4d+1)+K(23d^2)+(64d^3).
\]
This equation is quadratic in $K$ with discriminant $-64d^3+4d+1$ which is
negative for any integer $d>0$ which excludes the existence of
a power integral basis.

\noindent
{\bf I.2.} Let $-d\equiv 1 \ (\bmod \ 4), -d\neq -3$.
Taking any of $(Y_1,Y_2)=(1,-1),(1,0),(0,1)$ we have the following arguments.
$J(\vartheta)$ is divisible
by $y_0^3$, therefore $y_0=\pm 1$. The remaining factor $J_1$ of 
$J(\vartheta)$ is quadratic in $y_0$, hence we may substitute $y_0=1$
to get $J_2$. 
Then $J_2\pm 1=0$ must hold for a solution. 
Set $K=a^2+3a+9$ then
\[
J_2=(d+1)K^2+2d^2K+d^3.
\]
Equation $J_2-1=0$ has discriminant $-d^3+d+1$, 
equation $J_2+1=0$ has discriminant $-d^3-d-1$
in $K$. Both of them are negative for any $d>0$, $-d\equiv 1 \ (\bmod \ 4)$.

\vspace{0.5cm}
\noindent
{\bf II.} Consider the Gaussian integers $(Y_1,Y_2)=(i,-i),(i,0),(0,i)$.
In this case $d=1$. For each pair $(Y_1,Y_2)$ we calculate $X_1,X_2$ and take 
$\varepsilon=1$ and $\varepsilon=i$ in (\ref{ttee})
(these are all units up to sign). 
Then $J(\vartheta)$ depends only on $a$ and $y_0$. In each case
$J(\vartheta)$ has two factors, say $G_1,G_2$, one of being a complete square.
We solved the polynomial equations $G_1=\pm 1$, $G_2=\pm 1$ by using resolvents
and found the solutions listed in Theorem \ref{t2}.

\noindent
{\bf II.1.} Similarly we tested $(Y_1,Y_2)=(1,-1),(1,0),(0,1)$ with $-d=-1$
and $\varepsilon=1,i$ which did not give any solutions.

\vspace{0.5cm}
\noindent
{\bf III.} Consider the Eulerian integers $(Y_1,Y_2)=
(\omega_3,-\omega_3),(0,\omega_3),(\omega_3,0),\\(1-\omega_3,0),(0,1-\omega_3),
(1-\omega_3,-1+\omega_3)$. In this case $d=3$.
For each pair $(Y_1,Y_2)$ we calculate $X_1,X_2$ and take 
$\varepsilon=1,\omega_3,\omega_3^2$ in (\ref{ttee})
(these are all units up to sign).
Then $J(\vartheta)$ depends only on $a$ and $y_0$. In each case
$J(\vartheta)$ has two factors, say $G_1,G_2$.
Solving the polynomial equations $G_1=\pm 1$, $G_2=\pm 1$ by using resolvents
we did not find any integer solutions.

\noindent
{\bf III.1.} Similarly we tested $(Y_1,Y_2)=(1,-1),(1,0),(0,1)$ with $-d=-3$
and $\varepsilon=1,\omega_3,\omega_3^2$ which did not give any solutions.

\vspace{0.5cm}
\noindent
{\bf IV.} Consider the solutions $(a,Y_1,Y_2)$ given in Lemma \ref{lemma1}.
The parameter $a$ is given, we calculate $X_1=Y_1-aY_2$, $X_2=Y_2$. 
Using suitable units $\varepsilon$ we calculate $J(\vartheta)$
which always had two factors $G_1,G_2$. 
Solving the polynomial equations $G_1=\pm 1$, $G_2=\pm 1$ by using resolvents
we did not find any integer solutions.

\noindent
This proves Theorem \ref{t2}. $\Box$

\vspace{0.5cm}
\noindent
{\bf Computational remarks.} All our computations were performed in Maple \cite{maple}
and was executed on an average laptop. The calculations took some seconds only.

\end{document}